\documentclass[12pt, a4 paper]{amsart}
\usepackage{amssymb,amsmath,amsfonts,latexsym, tikz}
\usepackage{bm}
\usepackage{comment}
\usepackage{xcolor}
\usepackage[english]{babel}
\makeatletter\let\@secnumpunct\@empty\makeatother
\usepackage{titlesec}
\usepackage{url}
\newcommand{\bea}{\begin{eqnarray}}
\newcommand{\eea}{\end{eqnarray}}

\def\textmatrix#1&#2\\#3&#4\\{\bigl({#1 \atop #3}\ {#2 \atop #4}\bigr)}
\def\dispmatrix#1&#2\\#3&#4\\{\left({#1 \atop #3}\ {#2 \atop #4}\right)}
\newcommand{\be}{\begin{equation}}
\newcommand{\ee}{\end{equation}}
\newcommand{\ben}{\begin{eqnarray*}}
\newcommand{\een}{\end{eqnarray*}}

\newcommand{\bi}{\begin{itemize}}
\newcommand{\ei}{\end{itemize}}

\usepackage{titlesec}

\titleformat{\subsection}[block]
{\normalfont\large\bfseries}
{\thesubsection.}{1em}{}

\usepackage{amsthm}
\newtheorem{theorem}{Theorem}[section]
\newtheorem{lemma}[theorem]{Lemma}
\newtheorem{corollary}[theorem]{Corollary}

\usepackage{enumitem}
\usepackage{upgreek}
\usepackage{mathrsfs}

\usepackage{wasysym}
\usepackage[thinc]{esdiff}
\newtheorem{remark}{Remark}

\numberwithin{equation}{section}

\let\phi=\varphi

\begin{document}

\title[Algebras Generated by Two Idempotent]{Algebras Generated by Two Idempotent Matrices: Classification up to Isomorphism and Applications}

\author[Biswas]{Rounak Biswas}
\address{Department of Basic Science and Humanities, Indian Institute of Information Technology Bhagalpur, Bihar-813210, India.}
\email{xyzrounak3@gmail.com}



\begin{abstract}
This paper is a continuation of the earlier work by the author and Falguni Roy \cite{biswas2025comprehensive} on the algebra \(\operatorname{alg}(P,Q)\) generated by two idempotent matrices \(P,Q\in\mathbb{C}^{n\times n}\). A complete
classification of these algebras up to algebra isomorphism is accomplished, without imposing any
additional restrictions on the idempotents. The classification can be reduced to the case
of arbitrary pairs of idempotent matrices satisfying the condition that $\sigma\bigl((P-Q)^2\bigr)$
is a singleton. A minimal atlas of concrete algebras is presented such that every algebra
subject to this condition is isomorphic to exactly one member of this atlas. As an application, group invertibility in \(\operatorname{alg}(P,Q)\) is employed to derive an explicit solution to a particular system of linear equations.
\end{abstract}
\maketitle
\section{~Introduction}
Let \(P,Q\in\mathbb{C}^{n\times n}\) be idempotent matrices; that is,
\(P^2=P\) and \(Q^2=Q\). The present article is concerned with
the classification, up to algebra isomorphism, of the smallest
subalgebra of \(\mathbb{C}^{n\times n}\) containing both \(P\)
and \(Q\), which is denoted by \(\operatorname{alg}(P,Q)\).
A central motivation for investigating these algebras is to characterize generalized invertibility of their elements. 
In this direction, B{\"o}ttcher and Spitkovsky studied group invertibility in \(\operatorname{alg}(P,Q)\) for tightly coupled idempotent matrices and classified the corresponding algebras \cite{bottcher2011certain,bottcher2012group,bottcher2013classification}. Recall that an element \(A\in\mathcal A\) is called group invertible  if there exists \(A^g\in\mathcal A\) satisfying
\[
AA^g=A^gA,\qquad (A^g)^2A=A^g,\qquad A^2A^g=A~ \text{ (see \cite{kirkland2012group})}.
\] Subsequently, the author and Falguni Roy extended this study to non-tightly coupled idempotents satisfying \((PQ)^m=(PQ)^{m-1}\) for some \(m\in\mathbb N\) \cite{biswas2024drazin}. Recently, the author and Falguni Roy \cite{biswas2025comprehensive} obtained a complete classification of \(\operatorname{alg}(P,Q)\) for arbitrary pairs of idempotent matrices, without imposing any additional assumptions on \(P\) and \(Q\).  The notation and terminology used in the present article follow those adopted in \cite{biswas2025comprehensive}.
In particular, Lemma 2.2 of \cite{biswas2025comprehensive} establishes that if
\[
\sigma\bigl((P-Q)^2\bigr)=\{\lambda_1,\lambda_2,\ldots,\lambda_k\},
\]
then there exist \(k\) pairs of idempotent matrices \((P_1,Q_1),\ldots,(P_k,Q_k)\), of suitable orders, such that
\[
\sigma\bigl((P_i-Q_i)^2\bigr)=\{\lambda_i\},
\qquad i=1,\ldots,k,
\]
and
\[
P=\begin{bmatrix}
P_1&&\\
&\ddots&\\
&&P_k
\end{bmatrix},
\qquad
Q=\begin{bmatrix}
Q_1&&\\
&\ddots&\\
&&Q_k
\end{bmatrix}.
\]
Consequently,
\[
\begin{aligned}
\operatorname{alg}(P,Q)
=\operatorname{alg}(P_1,Q_1)
\oplus\operatorname{alg}(P_2,Q_2)\oplus\cdots\oplus\operatorname{alg}(P_k,Q_k).
\end{aligned}
\]
Thus, in order to study the algebra generated by an arbitrary pair of idempotent matrices, it is sufficient to consider the case in which \((P-Q)^2\) has a singleton spectrum.
For such a pair of idempotent matrices satisfying
\[
\sigma\bigl((P-Q)^2\bigr)=\{\lambda\},
\qquad \lambda\in\mathbb C,
\]
the possible algebraic structures of \(\operatorname{alg}(P,Q)\) have already been investigated in \cite{biswas2025comprehensive}. As noted in \cite{biswas2025comprehensive}, if
$P,Q\in\{0_n,I_n\}$, then $\operatorname{alg}(P,Q)$ is a
trivial algebra. Similarly, if $P=Q$, then
$\operatorname{alg}(P,Q)=\mathbb{C}P$. To exclude these
trivial cases, throughout this article, we restrict our
attention to pairs of distinct, nonzero, nonidentity
idempotent matrices; that is, $P\neq Q$  and 
$P,Q\notin\{0_n,I_n\}$.

In his review \cite{bottcher2025review} of \cite{biswas2025comprehensive}, Albrecht B{\"o}ttcher wrote,
``A complete classification of the algebras ${\rm alg}(P, Q)$ in the sense of providing a
minimal atlas of concrete algebras such that each of the algebras under consideration
with singleton spectrum of $(P - Q)^2$ is isomorphic to exactly one algebra of the atlas is
left as a subject of future research.''
The present article provides such an atlas. It builds upon the classification of \cite{biswas2025comprehensive} and presents a classification of the resulting algebras up
to isomorphism. 

The atlas consists of the $5$ infinite series of algebras (\ref{1})–(\ref{5}) listed in the following. 
For \(k\geq 1\), let \(\mathcal{B}_k\) be the algebra
\(
M_2\bigl(\mathbb{C}[t]/\langle t^k\rangle\bigr),
\)
consisting of all \(2\times 2\) matrices with entries in the quotient algebra \(\mathbb{C}[t]/\langle t^k\rangle\). Then every algebra \(\operatorname{alg}(P,Q)\) satisfying
\[
\sigma\bigl((P-Q)^2\bigr)=\{\lambda\},
\qquad \lambda\in\mathbb{C}\setminus\{0,1\},
\]
is isomorphic to exactly one algebra in the series
\begin{equation}
\mathcal{B}_1,\mathcal{B}_2,\mathcal{B}_3,\mathcal{B}_4,\ldots.
\label{1}
\end{equation}
Furthermore, we will show that the index \(k\) in \(\mathcal{B}_k\) coincides with the nilpotency index of \(\lambda I_n-(P-Q)^2\).
Next, in the case \(\sigma\bigl((P-Q)^2\bigr)=\{0\}\), three infinite series of algebras are constructed:
\begin{align}
& \mathcal{C}_2^{\mathrm{I}},
\mathcal{C}_3^{\mathrm{I}},\mathcal{C}_4^{\mathrm{I}},\ldots,
\label{2}\\
& \mathcal{C}_2^{\mathrm{II}},
\mathcal{C}_3^{\mathrm{II}},\mathcal{C}_4^{\mathrm{II}},\ldots,
\label{3}\\
& \mathcal{C}_2^{\mathrm{III}},
\mathcal{C}_3^{\mathrm{III}},\mathcal{C}_4^{\mathrm{III}},\ldots.
\label{4}
\end{align}
These series have the property that every algebra \(\operatorname{alg}(P,Q)\) for which \(P-Q\) is nilpotent of index \(k\) is isomorphic to exactly one of the three algebras \(\mathcal{C}_k^{\mathrm{I}},\mathcal{C}_k^{\mathrm{II}}\), and \(\mathcal{C}_k^{\mathrm{III}}\).
Finally, in the case \(\sigma\bigl((P-Q)^2\bigr)=\{1\}\), the classification of B{\"o}ttcher and Spitkovsky \cite{bottcher2013classification} is employed to establish that every algebra \(\operatorname{alg}(P,Q)\) is isomorphic to exactly one algebra in a series
\begin{align}
 \mathcal{A}_2,\mathcal{A}_3,\mathcal{A}_4,\ldots
\label{5}
\end{align} we will construct in Section \ref{2.3}.
It is also established that no two distinct algebras among the five series listed in (\ref{1})–(\ref{5}) are isomorphic.

As an application of the results on group invertibility in
$\operatorname{alg}(P,Q)$ established in
\cite{biswas2025comprehensive}, an explicit solution to the
system of linear equations \eqref{sole} is obtained by
employing the group inverse of a suitable element in the
algebra generated by two idempotent matrices associated
with the system. The system \((\ref{sole})\) was originally introduced by L.~M.~Hern\'andez-Ramos \cite{hernandez2005alternating}, who proposed an iterative procedure, known as the alternating oblique projections (AOP) method, for computing its solution. In contrast, the present approach yields an explicit expression for the solution in terms of the group inverse. We refer the reader to \cite{bottcher2010gentle,roch2010non} for a more comprehensive account of the properties and results involving pairs of idempotents.

The paper is organized as follows. Section \(\ref{s2}\) is devoted to the classification of \(\operatorname{alg}(P,Q)\) up to isomorphism. It is divided into three subsections, corresponding to the cases \(\sigma((P-Q)^2)=\{\lambda\}\), where \(\lambda\in\mathbb C\setminus\{0,1\}\), \(\sigma((P-Q)^2)=\{0\}\), and \(\sigma((P-Q)^2)=\{1\}\), respectively. Section \(\ref{s3}\) presents an application of group invertibility in \(\operatorname{alg}(P,Q)\)  to obtain an explicit solution of a particular system of linear equations.

\section{~Classification of the algebra up to isomorphism}\label{s2}
The aim of this section is to characterize, up to algebra isomorphism, \(\operatorname{alg}(P,Q)\), where \(P\) and \(Q\) are idempotent matrices satisfying
\[
\sigma\bigl((P-Q)^2\bigr)=\{\lambda\}.
\]
We divide the investigation into three cases: \(\lambda\in\mathbb C\setminus\{0,1\}\), \(\lambda=0\), and \(\lambda=1\). Throughout this section, for a given algebra \(\mathcal A\), we denote its center by \(C_{\mathcal A}\) and the set of all its nilpotent elements by \(N_{\mathcal A}\). Furthermore, by \(\mathbb C[t]/\langle(t-\lambda)^k\rangle\), we denote the quotient algebra of the polynomial algebra \(\mathbb C[t]\) by the ideal generated by \((t-\lambda)^k\), where \(\lambda\in\mathbb C\). Accordingly, \(M_2\bigl(\mathbb C[t]/\langle(t-\lambda)^k\rangle\bigr)\) denotes the algebra of all \(2\times2\) matrices whose entries belong to this quotient algebra. 
\subsection{$\sigma((P-Q)^2)=\{\lambda\}$ and $\lambda\in\mathbb{C}\setminus\{0,1\}$}
Throughout this section, we consider a pair of idempotent matrices \(P\) and \(Q\) such that \(\sigma\bigl((P-Q)^2\bigr)=\{\lambda\}\), where \(\lambda\in\mathbb{C}\setminus\{0,1\}\).
\begin{lemma}
    Let $P,Q\in \mathbb{C}^{n\times n}$ be two idempotent matrices with $\sigma\bigl((P-Q)^2\bigr)=\{\lambda\},$ and let the nilpotency index of $\lambda I_n-(P-Q)^2$ be $k$. If $\lambda\in\mathbb{C}\setminus\{0,1\}$,  then ${\rm alg}(P,Q)$ is a unital algebra.
    \begin{proof}
        The result follows from the fact that $\bigl(\lambda I_n-(P-Q)^2\bigr)^k=0_n$.
    \end{proof}
\end{lemma}
\begin{lemma}\label{center}
    Let $P,Q\in \mathbb{C}^{n\times n}$ be two idempotent matrices with $\sigma\bigl((P-Q)^2\bigr)=\{\lambda\},$ and let the nilpotency index of $\lambda I_n-(P-Q)^2$ be $k$. If $\lambda\in\mathbb{C}\setminus\{0,1\}$, then $C_{{\rm alg}(P,Q)}={\rm alg}\bigl(I_n,(P-Q)^2\bigr)$.
    \begin{proof}
        Let \(X=(P-Q)^2\). Then $XP=PX$ and $QX=XQ$. By employing the identities
\[
PQP=P(I_n-X),\qquad QPQ=Q(I_n-X),
\]
we deduce that every element \(W\in\operatorname{alg}(P,Q)\) can be expressed in the form
\[
W=a(X)+b(X)P+c(X)Q+d(X)PQ+e(X)QP,
\]
where \(a,b,c,d\), and \(e\) are polynomials. Now, suppose that \(W\in C_{\operatorname{alg}(P,Q)}\). Then, using the commutativity relation \(WP=PW\), we obtain \begin{align}\label{EQ1}
            (c+e)QP+(d-e)PQP-(c+d)PQ=0_n.
        \end{align} Multiplying (\ref{EQ1}) by \(P\) from the left, we obtain\[
-(c+d)PQP+(c+d)PQ=(c+d)PQ(I_n-P)=0_n.
\]
Consequently, multiplying the preceding equation on the left by \((I_n-P)Q\), we obtain
\[ 
(c+d)(I_n-P)QPQ(I_n-P)=0_n.
\]
Now, using the identity
\[
X(I_n-X)(I_n-P)=(I_n-P)QPQ(I_n-P),
\]
we obtain
\[
(c+d)X(I_n-X)(I_n-P)=0_n.
\]
Since \(\sigma(X)=\{\lambda\}\), with \(\lambda\notin\{0,1\}\), both \(X\) and \(I_n-X\) are invertible. Thus, the preceding equality implies that
\(
(c+d)(I_n-P)=0_n.
\) Furthermore, multiplying this equality by \(Q\) from both sides, we obtain
\[
(c+d)Q(I_n-P)Q=(c+d)QX=0_n.
\]
Since \(X\) is invertible and commutes with \(Q\), it follows that
\(
(c+d)Q=0_n.
\)
Consequently, using \(X=P+Q-PQ-QP\), we have \begin{align*}
            (c+d)X=&(c+d)(P+Q-PQ-QP)\\
            &=c+d,
        \end{align*} which further implies that
\(
(c+d)(X-I)=0_n.
\)
Since \(X-I\) is invertible, it follows that \(c+d=0\), or equivalently,
\(
c=-d.
\)
Substituting this relation into (\ref{EQ1}), we get
\[
(c+e)(QP-PQP)=0_n.
\]
Proceeding in a similar manner as above, one can show that this further implies
\(
c+e=0.
\) Hence we have \begin{align*}
            W&=a(X)+b(X)P+c(X)(Q-PQ-QP)\\
            &=a(X)+b(X)P+c(X)(X-P)\\
            &=a(X)+Xc(X)+P(b(X)-c(X))\\
            &=\alpha(X)+\beta(X)P,
        \end{align*} for some polynomials $\alpha$ and $\beta$. Finally, using the commutativity relation \(QW=WQ\), we obtain
\[
\beta(PQ-QP)=0_n.
\]
Since \(PQ-QP\) is invertible by \cite[Theorem 2.1]{biswas2025comprehensive}, we must have \(\beta=0_n\). Consequently, \(W=\alpha(X)\). Thus, every element in the center of \(\operatorname{alg}(P,Q)\) is a polynomial in \(X\). Hence, we conclude that
\[
C_{{\rm alg}(P,Q)}
=\operatorname{alg}(I_n,X)
=\operatorname{alg}\bigl(I_n,(P-Q)^2\bigr).
\] 
    \end{proof}

\end{lemma}

\begin{theorem}\label{lambdanonzeroalg}
    Let $P,Q\in \mathbb{C}^{n\times n}$ be two idempotent matrices such that $\sigma\bigl((P-Q)^2\bigr)=\{\lambda\},$ where $\lambda\in \mathbb{C}\setminus \{0,1\}$. If the nilpotency index of $\lambda I_n-(P-Q)^2$ is $k$,  then $${\rm alg}(P,Q)\cong M_2(\mathbb{C}[t]/\langle (t-\lambda)^k\rangle).$$ 
\begin{proof}
    Let \(X=(P-Q)^2\) and \(U=X(I_n-X)\). Since \(\sigma(X)=\{\lambda\}\), where \(\lambda\in\mathbb{C}\setminus\{0,1\}\), both \(X\) and \(I_n-X\) are invertible. Consequently, \(U\) is also invertible. We now define the following matrices
    \begin{align*}
        E_{11}&=P,~E_{22}=I_n-P,~E_{12}=PQ(I_n-P) \text{ and }\\ E_{21}&=U^{-1}(I_n-P)QP.
    \end{align*} Then \begin{align}\label{M1}
        E_{12}E_{21}&=PQ(I_n-P)U^{-1}(I_n-P)QP\notag\\
        &=U^{-1}(PQ-PQP)QP\notag\\
        &=U^{-1}(PQP-PQPQP)\notag\\
        &=U^{-1}(P(I_n-X)-P(I_n-X)^2)\notag\\
        &=PU^{-1}(X(I_n-X))=P=E_{11}. 
    \end{align} Similarly \begin{align}\label{M2}
        E_{21}E_{12}&=U^{-1}(I_n-P)QPQ(I_n-P)\notag\\
        &=U^{-1}(I_n-P)(I_n-X)Q(I_n-P)\notag\\
        &=U^{-1}(I_n-X)(I_n-P)Q(I_n-P)\notag\\
        &=U^{-1}(I_n-X)(Q-QP-PQ+PQP)\notag\\
        &=U^{-1}(I_n-X)(X-P+P(I_n-X))\notag\\
        &=U^{-1}(I_n-X)X(I_n-P)=I_n-P=E_{22}.
    \end{align} Furthermore, straightforward calculations show that the following relations hold: 
\begin{equation}\label{M3}
\left\{
\hspace{0.5cm}
\begin{gathered}
E_{11}E_{22}=E_{22}E_{11}=0_n,\\
E_{11}E_{12}=E_{12},\qquad
E_{11}E_{21}=0_n,\\
E_{22}E_{12}=0_n,\qquad
E_{22}E_{21}=E_{21},\\
E_{12}E_{11}=0_n,\qquad
E_{12}E_{22}=E_{12},\\
E_{21}E_{11}=E_{21},\qquad
E_{21}E_{22}=0_n.
\end{gathered}
\right.
\end{equation}
Thus, \(E_{11},E_{12},E_{21},E_{22}\) satisfy the defining relations of the standard matrix units in \(M_2(\mathbb C)\). Moreover, since \(I_n=E_{11}+E_{22}\), we have \begin{align}\label{M4}
        Q&=(E_{11}+E_{22})Q(E_{11}+E_{22})\notag\\
        &=E_{11}QE_{11}+E_{11}QE_{22}+E_{22}QE_{11}+E_{22}QE_{22}\notag\\
        &= (I_n-X)E_{11}+E_{12}+X(I_n-X)E_{21}+XE_{22}.
    \end{align} 
    Define
\[
\phi:M_2\left(\mathbb C[t]/\langle(t-\lambda)^k\rangle\right)
\longrightarrow \operatorname{alg}(P,Q)
\]
by
\[
\begin{aligned}
\phi\left(\begin{bmatrix}1&0\\0&0\end{bmatrix}\right)&=E_{11},
&\phi\left(\begin{bmatrix}0&1\\0&0\end{bmatrix}\right)&=E_{12},\\
\phi\left(\begin{bmatrix}0&0\\1&0\end{bmatrix}\right)&=E_{21},
&\phi\left(\begin{bmatrix}0&0\\0&1\end{bmatrix}\right)&=E_{22},
\end{aligned}
\]
and
\[
\phi\left(\begin{bmatrix}t&0\\0&t\end{bmatrix}\right)=X.
\]
By relations \((\ref{M1})\), \((\ref{M2})\), and \((\ref{M3})\), the map \(\phi\) is a well-defined algebra homomorphism. Furthermore, by \((\ref{M4})\), both \(Q\) and \(P=E_{11}\)  belong to the range of \(\phi\). Hence, \(\phi\) is surjective.
To establish injectivity, suppose that
\[
\phi\left(
\begin{bmatrix}
p_{11}(t)&p_{12}(t)\\
p_{21}(t)&p_{22}(t)
\end{bmatrix}
\right)=0_n,
\]
where \(p_{ij}\in\mathbb \mathbb{C}[t]/\langle (t-\lambda)^k\rangle\). By the definition of \(\phi\), we have
\begin{align}\label{M5}
  p_{11}(X)E_{11}+p_{12}(X)E_{12}
+p_{21}(X)E_{21}+p_{22}(X)E_{22}=0_n.  
\end{align}
Multiplying \((\ref{M5})\) by \(E_{11}\) on the left and \(E_{21}\) on the right, and using the relations in \((\ref{M1})\), \((\ref{M2})\) and \((\ref{M3})\), we obtain
\[
p_{12}(X)E_{11}=0_n.
\]
Similarly, multiplying \((\ref{M5})\) by \(E_{21}\) on the left and \(E_{22}\) on the right gives
\[
p_{12}(X)E_{22}=0_n.
\]
Since \(E_{11}+E_{22}=I\), it follows that \(p_{12}(X)=0_n\). Proceeding similarly for the remaining entries, we obtain \(p_{ij}(X)=0_n\) for all \(i,j\in\{1,2\}\). Thus, \(\phi\) is injective. Consequently, alg$(P,Q)\cong M_2(\mathbb{C}[t]/\langle (t-\lambda)^k\rangle)$.
\end{proof}
    
\end{theorem}
\begin{lemma}
    For $\lambda\in \mathbb{C}$ and $k\in\mathbb{N}$, $\mathbb{C}[t]/\langle (t-\lambda)^k\rangle\cong \mathbb{C}[t]/\langle t^k\rangle$.
    \begin{proof}
        The map $\psi:\mathbb{C}[t]/\langle (t-\lambda)^k\rangle\mapsto \mathbb{C}[t]/\langle t^k\rangle$, defined by $$\psi(f(t))=f(t+\lambda), ~\text{ for }f\in \mathbb{C}[t]/\langle (t-\lambda)^k\rangle $$ is a well-defined algebra isomorphism.
    \end{proof}
\end{lemma}
As a consequence of the preceding lemma, we obtain the following result.
\begin{corollary}\label{isomorphic algebra}
    For $\lambda_1, \lambda_2\in \mathbb{C}\setminus \{0,1\}$ and $k\in\mathbb{N}$, $M_2(\mathbb{C}[t]/\langle (t-\lambda_1)^k\rangle)\cong M_2(\mathbb{C}[t]/\langle (t-\lambda_2)^k\rangle)$.
\end{corollary}
\begin{remark}\label{R1}
 \normalfont   Let \((P,Q)\) and \((M,N)\) be two pairs of idempotent matrices such that
\[
\sigma\bigl((P-Q)^2\bigr)=\{\lambda\},
\qquad
\sigma\bigl((M-N)^2\bigr)=\{\gamma\},
\]
where \(\lambda,\gamma\in\mathbb C\setminus\{0,1\}\). Suppose, moreover, that both \(\lambda I_n-(P-Q)^2\) and \(\gamma I_n-(M-N)^2\) are nilpotent matrices of nilpotency index \(k\). It follows from Theorems \ref{lambdanonzeroalg} and Corollary \ref{isomorphic algebra} that
\[
\operatorname{alg}(P,Q)\cong\operatorname{alg}(M,N).
\]
Hence, for each \(k\in\mathbb N\), the algebra generated by a pair of idempotent matrices satisfying the above conditions is unique up to algebra isomorphism.  
\end{remark}
Let $\mathcal{B}_k$ denote the algebra
$M_2\bigl(\mathbb{C}[t]/\langle t^k\rangle\bigr)$ of all
$2\times 2$ matrices over the quotient algebra
$\mathbb{C}[t]/\langle t^k\rangle$. Then, by Remark
\ref{R1} and Theorem \ref{lambdanonzeroalg}, for
$\lambda\in\mathbb{C}\setminus\{0,1\}$ and
$\sigma\bigl((P-Q)^2\bigr)=\{\lambda\}$, the algebra
$\operatorname{alg}(P,Q)$ is isomorphic to $\mathcal{B}_k$
if and only if the nilpotency index of
$\lambda I_n-(P-Q)^2$ is $k$. This establishes the
completeness and uniqueness of the classification
presented in the list (\ref{1}).
Furthermore, by combining Theorem \ref{lambdanonzeroalg} with Lemma \ref{center}, we obtain the following result.

\begin{lemma}
    Let $P,Q\in\mathbb{C}^{n\times n}$ be idempotent matrices such that
    $\sigma\bigl((P-Q)^2\bigr)=\{\lambda\}$, where
    $\lambda\in\mathbb{C}\setminus\{0,1\}$. If the nilpotency index of
    $\lambda I_n-(P-Q)^2$ is $k$, then
    \[
        \operatorname{alg}(P,Q)
        \cong M_2\bigl(C_{\operatorname{alg}(P,Q)}\bigr).
    \]
    \begin{proof}
        By Lemma \ref{center}, the map
        \[
            \phi:\mathbb{C}[t]\longrightarrow
            C_{\operatorname{alg}(P,Q)},
            \qquad
            \phi(f(t))=f(X),
        \]
        where $X=(P-Q)^2$ and $f(t)\in\mathbb{C}[t]$, is a
        surjective algebra homomorphism with kernel
        $\langle(t-\lambda)^k\rangle$. Hence, by the
        first isomorphism theorem,
        \[
            \mathbb{C}[t]/\langle(t-\lambda)^k\rangle
            \cong C_{\operatorname{alg}(P,Q)}.
        \]
        The desired result now follows from Theorem
        \ref{lambdanonzeroalg}.
    \end{proof}
\end{lemma}

\subsection{$\sigma((P-Q)^2)=\{0\}$}
We next consider a pair of idempotent matrices $P,Q$ satisfying
$\sigma\bigl((P-Q)^2\bigr)=\{0\}$. Suppose that $P-Q$ is
nilpotent of index $k$. In \cite{biswas2025comprehensive},
the algebra $\operatorname{alg}(P,Q)$ was classified into
three distinct types, namely, Type I, Type II, and Type III,
according to the relationship between
$P(P-Q)^{k-1}$ and $(P-Q)^{k-1}$. These three types are
characterized by the following conditions:
\[
\begin{aligned}
&\text{Type I:}  \text{ when } P(P-Q)^{k-1}\neq 0_n
\text{ and } (P-Q)^{k-1}\neq P(P-Q)^{k-1};\\
&\text{Type II:} \text{ when } P(P-Q)^{k-1}\neq 0_n
\text{ but } (P-Q)^{k-1}=P(P-Q)^{k-1};\\
&\text{Type III:}\text{ when } P(P-Q)^{k-1}=0_n.
\end{aligned}
\]
\begin{lemma}\label{nonunital}
    Let $P,Q\in \mathbb{C}^{n\times n}$ be two idempotent matrices with $\sigma\bigl((P-Q)^2\bigr)=\{0\}$. Then  ${\rm alg}(P,Q)$ is not a unital algebra. 
    \begin{proof}
    Since \(P\) and \(Q\) are idempotent matrices satisfying \(\sigma\bigl((P-Q)^2\bigr)=\{0\}\), it follows from \cite[Remark 2]{biswas2025comprehensive} that \(I_n\notin\operatorname{alg}(P,Q)\). Now, suppose that \(P-Q\) is nilpotent of nilpotency index \(k\). We first consider the Type I algebra, which arises when $P(P-Q)^{k-1}\neq 0$ and $(P-Q)^{k-1}\neq P(P-Q)^{k-1}$. Let \(E\in\operatorname{alg}(P,Q)\). Then by \cite[Theorem 4.3]{biswas2025comprehensive}, \(E\) admits a representation of the form $$E=\sum_{i=0}^{k-1}\alpha_i(P(P-Q)^{i})+\sum_{i=0}^{k-2}\beta_i(Q(P-Q)^i).$$ Now if $E$ is the identity element in alg$(P,Q)$, then from $E(P-Q)=P-Q$ it follows that \begin{align*}
        P-Q=&\sum_{i=0}^{k-3}\alpha_i(P(P-Q)^{i+1})+\sum_{i=0}^{k-3}\beta_i(Q(P-Q)^{i+1})\\&+(\alpha_{k-2}+\beta_{k-2})P(P-Q)^{k-1}.
    \end{align*} This contradicts the fact that 
\[
\bigl\{P,Q,P(P-Q),Q(P-Q),\ldots,
P(P-Q)^{k-1}\bigr\}
\] forms a basis for \(\operatorname{alg}(P,Q)\) (see \cite[Theorem 4.3]{biswas2025comprehensive}). Hence, the Type I algebra is non-unital. A similar argument establishes that the algebras of Type II and Type III are also non-unital.
    \end{proof}
\end{lemma}
For two idempotent matrices \(P,Q\) satisfying \(\sigma\bigl((P-Q)^2\bigr)=\{0\}\), suppose that \(P-Q\) is a nilpotent matrix of nilpotency index \(k\), that is, \((P-Q)^k=0\) and \((P-Q)^{k-1}\neq 0\). Then it was established in \cite{biswas2025comprehensive} that the dimensions of the Type I, Type II, and Type III alg$(P,Q)$ are \(2k-1\), \(2k-2\), and \(2k-2\), respectively. Consequently, a Type I algebra cannot be isomorphic to either a Type II or a Type III algebra. We now investigate the possibility of distinguishing Type II and Type III algebras up to isomorphism.

 \begin{theorem}\label{for zero}
    Let $P,Q\in \mathbb{C}^{n\times n}$ be two idempotents of ${\rm Type ~II}$ and $M,N\in \mathbb{C}^{n\times n}$ be two idempotents of ${\rm Type ~III}$, where both $P-Q$ and $M-N$ have the nilpotency index $k$. Then
\[
\operatorname{alg}(P,Q)\not\cong
\operatorname{alg}(M,N).
\] 
\begin{proof}
   Since \(P\) and \(Q\) are idempotents of Type II, it follows from \cite[Theorem 4.6]{biswas2025comprehensive} that every element \(A\in\operatorname{alg}(P,Q)\) admits a representation of the form
\[
A=\sum_{i=0}^{k-2}\alpha_i P(P-Q)^i
+\sum_{i=0}^{k-2}\beta_i Q(P-Q)^i,
\]
where \(\alpha_i,\beta_i\in\mathbb{C}\). By \cite[Theorem 2.1]{biswas2025comprehensive}, we have \(\sigma(A)=\{0,\alpha_0+\beta_0\}\). Therefore, \(A\) is nilpotent if and only if \(\alpha_0+\beta_0=0\). Thus, the collection of all nilpotent elements of \(\operatorname{alg}(P,Q)\) is given by
\begin{align*}
    N_{{\rm alg}(P,Q)}=&\text{span}\{P-Q,P(P-Q),Q(P-Q),\ldots,\\
    &P(P-Q)^{k-2},Q(P-Q)^{k-2}\}.
\end{align*}
 Now, setting $D=P-Q$, we have $D^k=0$ and $PD^i-QD^i=D^{i+1}$. Consequently, the above expression can be rewritten as  $$N_{{\rm alg}(P,Q)}=\text{span}\{D,D^2,\ldots,D^{k-1},PD,\ldots,PD^{k-2}\}.$$ Furthermore, using the identities
$DP=D-PD+D^2$ and  $D^2P=PD^2$, it follows that $$(N_{{\rm alg}(P,Q)})^{k-1}=\text{span}\{D^{k-1}, PD^{k-1}\}=\text{span}\{D^{k-1}\}.$$ Now recall that, for a Type II algebra, $P(P-Q)^{k-1}\neq 0_n$. Therefore, it follows that $${\rm alg}(P,Q)(N_{{\rm alg}(P,Q)})^{k-1}\neq\{0_n\}$$ in the Type II case.   
Similarly for $M,N$, if $E=M-N$, then one can verify here $$\bigl(N_{{\rm alg}(M,N)}\bigr)^{k-1}=\text{span }\{E^{k-1}\}.$$ Moreover, here using $M(M-N)^{k-1}=0_n$ it follows that $${\rm alg}(M,N)\bigl(N_{{\rm alg}(M,N)}\bigr)^{k-1}=\{0_n\}.$$
Thus, these two types of algebras can be distinguished by an algebraic property that is preserved under isomorphisms. We therefore arrive at the desired result.

\end{proof}
\end{theorem}

\begin{remark}\label{R2}
 \normalfont   It follows from Lemma \ref{nonunital} that an algebra
$\operatorname{alg}(P,Q)$ satisfying
$\sigma\bigl((P-Q)^2\bigr)=\{0\}$ cannot be isomorphic to
an algebra $\operatorname{alg}(M,N)$ satisfying
$\sigma\bigl((M-N)^2\bigr)=\{\lambda\}$, where
$\lambda\in\mathbb{C}\setminus\{0,1\}$ and $M,N$ are
idempotent matrices. Furthermore, Theorem \ref{for zero}
establishes that the algebras of Types I, II, and III
generated by pairs of idempotent matrices satisfying
$\sigma\bigl((P-Q)^2\bigr)=\{0\}$ are pairwise
non-isomorphic.
\end{remark}

For each $k\in\mathbb N$, let $\mathcal{C}_k^l$ represent the
algebra $\operatorname{alg}(P,Q)$ of Type $l$, where $P-Q$
is nilpotent of index $k$ and
$l\in\{\mathrm{I},\mathrm{II},\mathrm{III}\}$.
For $k>1$, the existence of
$\mathcal{C}_k^{\mathrm{I}}$, $\mathcal{C}_k^{\mathrm{II}}$,
and $\mathcal{C}_k^{\mathrm{III}}$ follows from Examples
4.4, 4.7, and 4.10 of \cite{biswas2025comprehensive},
respectively. For \(k=1\), we have \(P-Q=0_n\), so \(P=Q\). As discussed in the introduction, this gives rise to a trivial case.
Moreover, by Remark
\ref{R2}, no two distinct algebras in the lists
(\ref{2})--(\ref{4}) are isomorphic. Consequently, these
lists account for all possible algebras of the form
$\operatorname{alg}(P,Q)$ under consideration.

Hence, the algebras listed in (\ref{2})--(\ref{4})
constitute a complete classification, up to algebra
isomorphism, of all algebras $\operatorname{alg}(P,Q)$
generated by two idempotent matrices $P$ and $Q$ such that
$P-Q$ is nilpotent. Furthermore, by Remark \ref{R2}, none
of the algebras in the lists (\ref{2})--(\ref{4}) is
isomorphic to any algebra in the list (\ref{1}).

\subsection{$\sigma((P-Q)^2)=\{1\}$}\label{2.3}

As observed in \cite{biswas2025comprehensive}, if \(P,Q\in \mathbb{C}^{n\times n}\) are two idempotent matrices satisfying
\(
\sigma\bigl((P-Q)^2\bigr)=\{1\},
\)
then \(PQ\) is nilpotent. Furthermore, whenever \(PQ\) is nilpotent, let \(m\) be the least positive integer such that
\(
(PQ)^m=(QP)^m=0.
\)
A complete classification of \(\operatorname{alg}(P,Q)\) into four distinct types, based on the value of \(m\), was established in \cite{bottcher2013classification}.

\begin{lemma}\label{PQ1Unital}
    Let \(P,Q\in \mathbb{C}^{n\times n}\) be two idempotent matrices such that
$\sigma\bigl((P-Q)^2\bigr)=\{1\}$.
Then \(\operatorname{alg}(P,Q)\) is a unital algebra.
    \begin{proof}
        If $\sigma((P-Q)^2)=\{1\}$, then $(I_n-(P-Q)^2)^k=0$ for some integer $k$,  and the result follows.
    \end{proof}
\end{lemma}
\begin{remark}\label{R3}
  \normalfont  As a consequence of Lemma \ref{PQ1Unital} and Lemma \ref{nonunital}, it follows that if \((P,Q)\) and \((M,N)\) are pairs of idempotent matrices satisfying
\[
\sigma\bigl((P-Q)^2\bigr)=\{0\}
\quad\text{and}\quad
\sigma\bigl((M-N)^2\bigr)=\{1\},
\]
respectively, then
\[
\operatorname{alg}(P,Q)\not\cong
\operatorname{alg}(M,N).
\]
\end{remark}

\begin{theorem}\label{lambda,1}
   Let \((P,Q)\) and \((M,N)\) be two pairs of idempotent matrices satisfying
\(
\sigma\bigl((P-Q)^2\bigr)=\{1\}
\quad\text{and}\quad
\sigma\bigl((M-N)^2\bigr)=\{\lambda\},
\)
respectively, where \(\lambda\in\mathbb C\setminus\{0,1\}\). Then
\[
\operatorname{alg}(P,Q)\not\cong
\operatorname{alg}(M,N).
\]
    \begin{proof}
        Since \(P,Q\) are two idempotent matrices satisfying \(\sigma((P-Q)^2)=\{1\}\), it follows from \cite[Theorem 2.1]{biswas2025comprehensive} that there exist two continuous algebra homomorphisms $G_1$ and $G_2$ from \(\operatorname{alg}(P,Q)\) onto \(\mathbb {C}\). Now if there exists a continuous algebra homomorphism $$\phi:{\rm alg}(M,N)\mapsto \mathbb{C},$$ then $\phi(M),\phi(N)\in \{0,1\}$. This further implies that $\phi((M-N)^2)\in\{0,1\}$. On the other hand, $\lambda I_n-(M-N)^2$ is a nilpotent matrix; thus $\phi(\lambda I_n-(M-N)^2)=0$, which contradicts that $\lambda\in \mathbb{C}\setminus\{0,1\}$. Hence, the desired result follows.
    \end{proof}
\end{theorem}
Finally, let \(M,N\in\mathbb{C}^{n\times n}\) be idempotent matrices satisfying
\[
(MN)^m=(NM)^m=0_n
\]
for some \(m\in\mathbb{N}\). In \cite{bottcher2013classification}, B{\"o}ttcher and Spitkovsky classified the algebra \(\operatorname{alg}(M,N)\) into four distinct types, up to algebra isomorphism. They constructed concrete algebras $Z_k$
such that \(\operatorname{alg}(M,N)\) is isomorphic to exactly one of the algebras in the  list
\begin{align}\label{zero}
    Z_0,Z_1,Z_2,Z_3,Z_4,\ldots.
\end{align}
Now, let \(P,Q\in\mathbb{C}^{n\times n}\) be idempotent matrices such that
\[
\sigma\bigl((P-Q)^2\bigr)=\{1\}.
\]
Then there exists a least \(k\in\mathbb{N}\) such that
\[
(PQ)^k=(QP)^k=0_n.
\]
Consequently, \(\operatorname{alg}(P,Q)\) is isomorphic to exactly one of the algebras in the list (\ref{zero}). However, \(\operatorname{alg}(P,Q)\) cannot be isomorphic to either \(Z_0\) or \(Z_1\), since
\[
\dim\operatorname{alg}(P,Q)>1
\]
whenever \(P,Q\notin\{0_n,I_n\}\) and \(P\neq Q\), whereas \(\dim Z_0=0\) and \(\dim Z_1=1\). Consequently, \(\operatorname{alg}(P,Q)\) is isomorphic to exactly one of the algebras
\[
Z_2,Z_3,Z_4,\ldots.
\]
For each \(i> 1\), let \(\mathcal{A}_i\) denote an algebra \(\operatorname{alg}(P,Q)\) satisfying \(\sigma\bigl((P-Q)^2\bigr)=\{1\}\) and \(\operatorname{alg}(P,Q)\cong Z_i\). By Theorem 2.4 of \cite{bottcher2013classification}, any two distinct algebras in the list (\ref{5}) are non-isomorphic. Finally, Remark \(\ref{R3}\) and Theorem \(\ref{lambda,1}\) ensure that no algebra in the list (\ref{5}) is isomorphic to any algebra in the lists (\ref{1})–(\ref{4}).
\section{~Application}\label{s3}
We now consider the following system of linear equations: \begin{align}\label{sole}
    \begin{bmatrix}
        A &B^*\\
        B&0_m
    \end{bmatrix}\begin{bmatrix}
        x\\y
    \end{bmatrix}=\begin{bmatrix}
        f\\0
    \end{bmatrix},
\end{align}
where $A\in \mathbb{C}^{n\times n}$ is a positive definite matrix and $B$ is a matrix of order $m\times n$ with rank $m$. In this section, $\ker (M)$, Im $(M)$ denote the kernel and the range of the matrix $M\in \mathbb{C}^{n_1\times n_2}$, respectively. By  \cite[Proposition 1.1]{hernandez2005alternating}, the system \eqref{sole} admits a unique solution  $(x^\prime,y^\prime)$ which satisfies $$Ax^\prime-f\in(\ker B)^\perp\text{ and }x^\prime\in \ker B.$$ Moreover, if we define $V:=\{x:f-Ax\perp \ker B\}$, then $\{x^\prime\}=V\cap \ker B.$ Since $A$ is a positive definite matrix, $\langle x,y\rangle_A:= \langle Ax,y\rangle=y^*Ax$ defines an inner product on $\mathbb{C}^n$. Furthermore, considering the orthogonal complement of \(\ker B\) with respect to this inner product, denoted by \((\ker B)^{\perp_A}\), it follows from \cite{hernandez2005alternating} that
\begin{align}\label{V}
    V=\{A^{-1}f\}+(\ker B)^{\perp_A}.
\end{align}

\begin{lemma}
    If $R$ is the orthogonal projection onto $(\ker B)^\perp$, and $Q=A^{-1}RA$, then \begin{enumerate}[label=\normalfont(\roman*)]
        \item $\ker (QR)=\ker (R)=\ker (B),$
        \item  ${\rm Im} ~(QR)={\rm Im}~(Q)=(\ker B)^{\perp_A}$.
    \end{enumerate} 
    \begin{proof}
        \begin{enumerate}[label=\normalfont(\roman*)]
            \item Clearly, \(\ker(R)\subseteq\ker(QR)\). To establish the reverse inclusion, let \(z\in\ker(QR)\). Then
\(
QRz=A^{-1}RARz=0.
\)
Consequently, \(RARz=0\), which implies
\(
\langle ARz,Rz\rangle=0.
\)
Since \(A\) is positive definite, we must have \(Rz=0\), and hence \(z\in\ker(R)\). Therefore, $\ker (QR)=\ker (R)$.
            \item Clearly, \(\operatorname{Im}(QR)\subseteq\operatorname{Im}(Q)\). Moreover, using the equality \(\ker(QR)=\ker(R)\), we obtain
\[
\begin{aligned}
\operatorname{rank}(QR)
&=n-\operatorname{nullity}(QR)\\
&=n-\operatorname{nullity}(R)\\
&=\operatorname{rank}(R).
\end{aligned}
\]
Now, since \(Q=A^{-1}RA\), we have \(\operatorname{rank}(Q)=\operatorname{rank}(R)\). Consequently,
$\operatorname{rank}(QR)=\operatorname{rank}(Q).
$
Together with \(\operatorname{Im}(QR)\subseteq\operatorname{Im}(Q)\), this gives
\(
\operatorname{Im}(QR)=\operatorname{Im}(Q).
\)

Next, let \(y\in\operatorname{Im}(Q)\). Then there exists \(x\) such that
\[
y=Qx=A^{-1}RAx.
\]

For any \(z\in\ker(B)\), we have
\[
\begin{aligned}
\langle y,z\rangle_A&= \langle Ay,z\rangle\\
&=\langle RAx,z\rangle\\
&=\langle Ax,Rz\rangle\\
&=0.
\end{aligned}
\]

Thus, \(\operatorname{Im}(Q)\subseteq(\ker B)^{\perp_A}\).

Conversely, let \(z\in(\ker B)^{\perp_A}\). Then \(Az\in(\ker B)^\perp\), and hence \(RAz=Az\). It follows that
\[
Qz=A^{-1}RAz=z,
\]

which implies \(z\in\operatorname{Im}(Q)\). Therefore, $\text{Im}(Q)=(\ker B)^{\perp_A}$.
\end{enumerate}
    \end{proof}
\end{lemma} 
 By the preceding lemma, we have
\begin{align*}
    \text{Im}(QR)\cap\ker (QR)&=\ker B\cap (\ker B)^{\perp_A}\\
    &=\{0\}.
\end{align*}
Consequently, \(QR\) is group invertible. Since \(Q\) and \(R\) are both idempotent matrices, the group inverse of $QR$ can be explicitly
computed using  \cite[Theorem 4.9]{biswas2025comprehensive}.
\begin{theorem}
    Let $(a,b)$ be the unique solution of $(\ref{sole})$. Then $a=A^{-1}f-(QR)^g(QR)A^{-1}f$ and $b=(BB^*)^{-1}B\bigl(A(QR)^gQRA^{-1}f\bigr)$.
    \begin{proof}
        If $(a,b)$ is the solution of (\ref{sole}), then $a\in V\cap \ker B$. Therefore, by (\ref{V}),  it follows that $a-A^{-1}f\in (\ker B)^{\perp_A}$.  Since $QR$ is group invertible, the matrix  $(QR)^gQR=QR(QR)^g$ is an idempotent matrix whose range coincides with  $\text{Im}(QR)=(\ker B)^{\perp_A}$. Hence we have \begin{align*}
            A^{-1}f-a&=(QR)^gQR(A^{-1}f-a)\\
            &=(QR)^gQR (A^{-1}f),      \end{align*}where the last equality follows from \(a\in\ker B=\ker(QR)\). Consequently, $a=A^{-1}f-(QR)^g(QR)A^{-1}f$. Finally, $b=(BB^*)^{-1}B\bigl(A(QR)^gQRA^{-1}f\bigr)$ follows from the fact that $B^* b=f-Aa$ and $BB^*$ is an invertible matrix.
    \end{proof}
\end{theorem}

\smallskip
\noindent\section*{Declaration of AI usage} ChatGPT Go was used solely to improve the language and presentation of the manuscript. All the results presented in this article were obtained without the use of any AI tools, and the author takes full responsibility for their correctness.
\smallskip

\section*{Declaration of competing interest}
There is no competing interest.
\smallskip
\section*{Data availability}
No data was used for the research described in the article.

\section*{Acknowledgement}
The author is grateful to Prof. Albrecht B{\"o}ttcher for his valuable comments and suggestions throughout the preparation of this article, which have contributed significantly to improving the manuscript.

\end{document}